\documentclass[11pt]{article}
\usepackage{amsmath}

\newtheorem{theorem}{Theorem}[section]

\newcommand\beq{\begin{equation}}
\newcommand\eeq{\end{equation}}
\newcommand\bce{\begin{center}}
\newcommand\ece{\end{center}}
\newcommand\bea{\begin{eqnarray}}
\newcommand\eea{\end{eqnarray}}
\newcommand\bean{\begin{eqnarray*}}
\newcommand\eean{\end{eqnarray*}}
\newcommand\ben{\begin{enumerate}}
\newcommand\een{\end{enumerate}}
\newcommand\bit{\begin{itemize}}
\newcommand\eit{\end{itemize}}
\newcommand\brr{\begin{array}}
\newcommand\err{\end{array}}
\newcommand\bt{\begin{tabular}}
\newcommand\et{\end{tabular}}

\newcommand\ms{\medskip}
\newcommand\ul{\underline}

\renewcommand\S{\mathcal S}
\newcommand\D{\mathcal D}
\newcommand\E{\mathcal E}
\newcommand\M{\mathcal M}
\newcommand\wh{\widehat}
\newcommand\wt{\widetilde}

\author{Emeric Deutsch~\thanks{Polytechnic Institute of New York University, Brooklyn, NY 11201.}
\and
Sergi Elizalde~\thanks{Department of Mathematics,
Dartmouth College, Hanover, NH 03755; sergi.elizalde@dartmouth.edu.}}
\title{The largest and the smallest fixed points of permutations}
\date{}

\begin{document}
\maketitle

\begin{abstract}
We give a new interpretation of the derangement numbers $d_n$ as the
sum of the values of the largest fixed points of all non-derangements of length $n-1$.
We also show that the analogous sum for the smallest fixed points equals the number of permutations of length $n$ with at least two fixed points.
We provide analytic and bijective proofs of both results, as well as a new recurrence for the derangement numbers.

\end{abstract}

\section{Largest fixed point}

Let $[n]=\{1,2,\dots,n\}$, and let $\S_n$ denote the set of permutations of $[n]$.
Throughout the paper, we will represent permutations using cycle notation unless specifically stated otherwise.
Recall that $i$ is a fixed point of $\pi\in\S_n$ if $\pi(i)=i$.
Denote by $\D_n$ the set of derangements of $[n]$, i.e., permutations with no fixed points, and let $d_n=|\D_n|$.
Given $\pi\in\S_n\setminus\D_n$, let $\ell(\pi)$ denote the largest fixed point of $\pi$.
Let $$a_{n,k}=|\{\pi\in\S_n : \ell(\pi)=k\}|.$$

Clearly, \beq\label{eq:ann} a_{n,1}=d_{n-1} \quad\mbox{and}\quad a_{n,n}=(n-1)!.\eeq It also follows from the definition that
\begin{equation}\label{eq:ank_sum}
a_{n,k}=d_{n-1}+ \sum_{j=1}^{k-1} a_{n-1,j},
\end{equation}
since by removing the largest fixed point $k$ of a permutation in $\S_n\setminus\D_n$, we get a permutation of $\{1,\dots,k-1,k+1,\dots,n\}$
whose largest fixed point (if any) is less than $k$.
If in (\ref{eq:ank_sum}) we replace $k$ by $k-1$, then by subtraction we obtain
\begin{equation}\label{eq:ank_rec}
a_{n,k}=a_{n,k-1}+a_{n-1,k-1}
\end{equation}
for $k\ge2$, or equivalently, $a_{n,k}=a_{n,k+1}-a_{n-1,k}$ for $k\ge1$.
Together with the second equation in~(\ref{eq:ann}), it follows that the numbers
$a_{n,k}$ form Euler's difference table of the factorials (see~\cite{DuRa,Ges,Rak}).
Table~\ref{tab:ank} shows the values of $a_{n,k}$ for small $n$.
The combinatorial interpretation given in~\cite{DuRa,Ges} is that $a(n,k)$ is
the number of permutations of $[n-1]$ where none of $k,k+1,\dots,n-1$ is a fixed point.
This interpretation is clearly equivalent to ours using the same reasoning behind equation~(\ref{eq:ank_sum}).

\begin{table}[htb]
$$ \begin{array}{c|rrrrrr}
n\backslash k & 1 & 2 & 3 &4 &5 & 6\\ \hline
1 & 1   \\
2 & 0 & 1 \\
3 & 1 & 1 & 2 \\
4 & 2 & 3 & 4 & 6 \\
5 & 9 & 11 & 14 & 18 & 24\\
6 & 44 & 53 & 64 & 78 & 96 & 120
\end{array}$$
\caption{\label{tab:ank} The values of $a_{n,k}$ for $n$ up to 6.}
\end{table}

We point out that it is possible to give a direct combinatorial proof of the recurrence~(\ref{eq:ank_rec})
from our definition of the $a_{n,k}$.
Indeed, let $\pi\in\S_n$ with $\ell(\pi)=k$. If $\pi(1)=m\neq1$, then the permutation of $[n]$ obtained from the one-line notation of $\pi$
by moving $m$ to the end, replacing $1$ with $n+1$, and subtracting one from all the entries has largest fixed point $k-1$.
If $\pi(1)=1$, then removing $1$ and subtracting one from the remaining entries of $\pi$ we get
a permutation of $[n-1]$ whose largest fixed point is $k-1$.

\ms

Define
\begin{equation}\label{eq:alpha}
\alpha_n=\sum_{k=1}^n k a_{n,k}=\sum_{\pi\in\S_n\setminus\D_n} \ell(\pi).
\end{equation}

We now state our main result, which we prove analytically and bijectively in the next two subsections.
\begin{theorem}\label{thm:main}
For $n\ge1$, we have $$\alpha_n=d_{n+1}.$$
\end{theorem}

\subsection{Analytic proof}
Replacing $n$ by $n+1$, from (\ref{eq:alpha}) we have
\begin{equation}\label{eq:alpha1}
\alpha_{n+1} = a_{n+1,1}+2a_{n+1,2}+\dots+na_{n+1,n}+(n+1)a_{n+1,n+1}.
\end{equation}
Adding (\ref{eq:alpha}) and (\ref{eq:alpha1}) and taking into account (\ref{eq:ank_rec}), we obtain
\begin{equation}\label{eq:alpha2}
\alpha_n + \alpha_{n+1} = a_{n+1,2}+2a_{n+1,3}+\dots+na_{n+1,n+1}+(n+1)!.
\end{equation}
Adding (\ref{eq:alpha2}) with the obvious equality
$$(n+1)! - d_{n+1}=a_{n+1,1}+a_{n+1,2}+\dots+a_{n+1,n}+a_{n+1,n+1},$$
we obtain
$$\alpha_n + \alpha_{n+1} +(n+1)! - d_{n+1} = \alpha_{n+1} + (n+1)!,$$
whence $\alpha_n = d_{n+1}$.

\subsection{Bijective proof}

To find a bijective proof of Theorem~\ref{thm:main}, we first construct a set whose cardinality is $\alpha_n$.
Let $\M_n\subset (\S_n\setminus\D_n)\times [n]$ be the set of pairs $(\pi,i)$ where $\pi\in\S_n\setminus\D_n$ and $i\le\ell(\pi)$.
We underline the number $i$ in $\pi$ to indicate that it is marked.
For example, we write $(2)(3)(7)(8)(1,\ul{4},9)(5,6)$ instead of the pair $((2)(3)(7)(8)(1,4,9)(5,6),4)$.
 It is clear that
$$|\M_n|=\sum_{k=1}^n k a_{n,k}=\alpha_n.$$
To prove Theorem~\ref{thm:main}, we give a bijection between $\D_{n+1}$ and $\M_n$.

Given $\pi\in\D_{n+1}$, we assign to it an element $\wh\pi\in\M_n$ as follows. Write $\pi$ as a product of cycles, starting with
the one containing $n+1$, say $$\pi=(n+1,i_1,i_2,\dots,i_r)\,\sigma.$$
Let $q$ be the largest index, $1\le q\le r$, such that $i_1<i_2<\dots<i_q$.
We define
$$\wh\pi=\begin{cases} (i_1)(i_2)\dots(\ul{i_r})\,\sigma & \mbox{if } q=r,\\
(i_1)(i_2)\dots(i_q)(\ul{i_{q+1}},i_{q+2},\dots,i_r)\,\sigma & \mbox{if } q<r.\end{cases}$$

Now we describe the inverse map. Given $\wh\pi\in\M_n$, let its unmarked fixed points be $i_1<i_2<\dots<i_q$, and let $j_1$ be the marked element.
We can write $\wh\pi=(i_1)\dots(i_q)(\ul{j_1},j_2,\dots,j_t)\,\sigma$. Notice that $t=1$ if the marked element is a fixed point. Define
$$\pi=(n+1,i_1,i_2,\dots,i_q,j_1,j_2,\dots,j_t)\,\sigma.$$

Here are some examples of the bijection between $\D_{n+1}$ and $\M_n$:
\bean\pi=(12,2,4,9,7,5,6)(1,3)(8,11,10) &\leftrightarrow& \wh\pi=(2)(4)(9)(\ul{7},5,6)(1,3)(8,11,10),\\
\pi=(10,2,7,8,3)(1,4,9)(5,6) &\leftrightarrow& \wh\pi=(2)(7)(8)(\ul{3})(1,4,9)(5,6),\\
\pi=(10,2,3,7,8,4,9,1)(5,6) &\leftrightarrow& \wh\pi=(2)(3)(7)(8)(\ul{4},9,1)(5,6). \eean

\section{Smallest fixed point}

In a symmetric fashion to the statistic $\ell(\pi)$, we can define
$s(\pi)$ to be the smallest fixed point of $\pi\in\S_n\setminus\D_n$.
Let $$b_{n,k}=|\{\pi\in\S_n : s(\pi)=k\}|.$$
The numbers $b_{n,k}$ appear in~\cite[pp. 174-176,185]{Cha} as $R_{n,k}$ (called {\em rank}).
Define
\begin{equation}\label{eq:beta}
\beta_n=\sum_{k=1}^n k b_{n,k}=\sum_{\pi\in\S_n\setminus\D_n} s(\pi).
\end{equation}

It is not hard to see by symmetry that
\beq\label{eq:ab}b_{n,k}=a_{n,n+1-k}.\eeq Indeed, one can use
the involution $\pi\mapsto\pi'$ on $S_n$ where $\pi'(i)=n+1-\pi(n+1-i)$.
Alternatively, another involution that proves~(\ref{eq:ab}) consists of replacing each entry $i$
in the cycle representation of $\pi\in\S_n$ by $n+1-i$; for example,
$(183)(2)(4975)(6)$ is mapped to $(927)(8)(6135)(4)$.

To find a combinatorial interpretation of $\beta_n$, let $\E_{n+1}$ be the set of permutations of $[n+1]$
that have at least two fixed points. We have that
\beq\label{eq:en}|\E_{n+1}|=(n+1)!-d_{n+1}-(n+1)d_n,\eeq
since out of the $(n+1)!$ permutations of $[n+1]$, there are $d_{n+1}$ derangements and $(n+1)d_n$ permutations having exactly one fixed point.

The following result is the analogue of Theorem~\ref{thm:main} for the statistic $s(\pi)$. We give an analytic proof based on
that theorem, and a directive bijective proof as well.
\begin{theorem}
For $n\ge1$, we have $$\beta_n=|\E_{n+1}|.$$
\end{theorem}

\subsection{Analytic proof}

From the definitions of $\alpha_n$ and $\beta_n$, and equation~(\ref{eq:ab}), it follows that
$$\alpha_n + \beta_n = (n+1)\sum_{k=1}^n a_{n,k}= (n+1)(n! - d_n).$$
Using Theorem~\ref{thm:main}, we have
$$\beta_n = (n+1)! - (n+1)d_n - d_{n+1},$$
which by~(\ref{eq:en}) is just the cardinality of $\E_{n+1}$ as claimed.

Note also the following identities involving $\beta_n$ which follow from the known recurrence $d_n=nd_{n-1}+(-1)^n$:
\bean &&\beta_n = (n+1)! + (-1)^n -2(n+1)d_n,\\
&&\beta_n = (n+1) \beta_{n-1} + n(-1)^{n+1}.\eean
The sequence $\beta_n$ starts $0,1,1,7,31,191,\dots$.
Using the well known fact that \beq\label{eq:limdn}\lim_{n\rightarrow\infty}\frac{d_n}{n!}=\frac{1}{e},\eeq
we see that $$\lim_{n\rightarrow\infty}\frac{\beta_n}{(n+1)!}=1 -\frac{2}{e}.$$

\subsection{Bijective proof}

Let $\M'_n\subset (\S_n\setminus\D_n)\times [n]$ be the set of pairs $(\pi,i)$ where $\pi\in\S_n\setminus\D_n$ and $i\le s(\pi)$.
As before, we underline the number $i$ in $\pi$ to indicate that it is marked.
It is clear that
$$|\M'_n|=\sum_{k=1}^n k b_{n,k}.$$

We now give a bijection between $\E_{n+1}$ and $\M'_n$.
Given $\pi\in\E_{n+1}$, let $i$ be its smallest fixed point. We can write
$$\pi=(i)(n+1,j_2,\dots,j_t)\,\sigma,$$
where no $j$s appear if $n+1$ is a fixed point. Define
$$\wt\pi=(\ul{i},j_2,\dots,j_t)\,\sigma.$$
Note that $\wt\pi\in\M'_n$, because if $\sigma$ has fixed points then they are all larger than $i$, and if it does not, then $t=1$ and $i$ is the smallest
fixed point of $\wt\pi$. Essentially, $\pi$ and $\wt\pi$ are related by conjugation by the transposition $(i,n+1)$.

Conversely, given $\wt\pi\in\M'_n$, let $i$ be the marked entry. We can write $$\wt\pi=(\ul{i},j_2,\dots,j_t)\,\sigma,$$
where no $j$s appear if $i$ is a fixed point. Then
$$\pi=(i)(n+1,j_2,\dots,j_t)\,\sigma.$$
Roughly speaking, we replace $\ul{i}$ with $n+1$ and add $i$ as a fixed point.
Note that if $t\ge2$ then $\sigma$ must have fixed points.

Here are some examples of the bijection between $\E_{n+1}$ and $\M_n$:
\bean
\pi=(3)(10,1,7,2,8)(5)(6)(4,9) &\leftrightarrow& \wt\pi=(\ul{3},1,7,2,8)(5)(6)(4,9),\\
\pi=(5)(10)(6)(3,1,7,2,8)(4,9) &\leftrightarrow& \wt\pi=(\ul{5})(6)(3,1,7,2,8)(4,9). \eean

\section{Other remarks}

\subsection{A recurrence for the derangement numbers}
An argument similar to the bijective proof of Theorem~\ref{thm:main} can be used to prove the recurrence \beq\label{eq:recdn}d_n=\sum_{j=2}^n(j-1)\binom{n}{j}d_{n-j}\eeq combinatorially as follows.

A derangement  $\pi\in\D_n$ can be written as a product of cycles, starting with the one containing $n$, say $$\pi=(n,i_1,i_2,\dots,i_r)\,\sigma.$$
Consider two cases: \bit
\item If $i_1<i_2<\dots<i_{r-1}$ (this is vacuously true for $r=1,2$),
then the number of choices for the numbers $i_1,\dots,i_r$ satisfying this condition is
$r\binom{n-1}{r}$, since we can first choose an $r$-subset of $[n-1]$ and then decide which one is $i_r$. Now, the number of choices
for $\sigma$ is $d_{n-r-1}$.
\item Otherwise, there is an index $1\le q\le r-1$ such that $i_1<i_2<\dots<i_q>i_{q+1}$.
In this case, there are $q\binom{n-1}{q+1}$ choices for the numbers $i_1,\dots,i_{q+1}$, since we can first
choose a $(q+1)$-subset of $[n-1]$ and then decide which element other than the maximum is $i_{q+1}$.
Now, there are $d_{n-q-1}$ choices for $(i_{q+1},\dots,i_r)\,\sigma$.
\eit
The total number of choices is
\begin{multline*}\sum_{r=1}^{n-1}r\binom{n-1}{r}d_{n-r-1}+\sum_{q=1}^{n-1}q\binom{n-1}{q+1}d_{n-q-1}
=\sum_{r=1}^{n-1}r\left(\binom{n-1}{r}+\binom{n-1}{r+1}\right)d_{n-r-1}\\
=\sum_{r=1}^{n-1}r\binom{n}{r+1}d_{n-r-1},\end{multline*}
which equals the right hand side of (\ref{eq:recdn}).

\ms

Alternatively, the recurrence~(\ref{eq:recdn}) is relatively straightforward to prove using generating functions.
Indeed, let $$D(x)=\sum_{n\ge0}d_n\frac{x^n}{n!}=\frac{e^{-x}}{1-x}$$ be the generating function for the number
of derangements. The generating function for the right hand side of~(\ref{eq:recdn}), starting from $n=1$, is
\begin{multline*}\sum_{n\ge1} \sum_{j=2}^n(j-1)\binom{n}{j}d_{n-j}\,\frac{x^n}{n!}
=\left(\sum_{i\ge0}d_i\frac{x^i}{i!}\right)\left(\sum_{j\ge1}(j-1)\frac{x^j}{j!}\right)\\
=\frac{e^{-x}}{1-x}\,(xe^x-e^x+1)=-1+\frac{e^{-x}}{1-x}=D(x)-1.\end{multline*}

\subsection{Probabilistic interpretation}
Let $X_n$ be the random variable that gives the value of the largest fixed point of a random element of $\S_n\setminus\D_n$.
Its expected value is then
$$E[X_n]=\frac{\sum_{k=1}^n k a_{n,k}}{|\S_n\setminus\D_n|}.$$ Theorem~\ref{thm:main} is equivalent to the fact that
\beq\label{eq:avg}E[X_n]=\frac{d_{n+1}}{n!-d_n}.\eeq
Using~(\ref{eq:limdn}), we get from equation~(\ref{eq:avg}) that
\beq\label{eq:limit}\lim_{n\rightarrow\infty} \frac{E[X_n]}{n}=\frac{1}{e-1}.\eeq

Occurrences of fixed points in a random permutation of $[n]$, normalized by dividing by $n$,
approach a Poisson process in the interval $[0,1]$ with mean 1 as $n$ goes to infinity.
An interpretation of equation~(\ref{eq:limit}) is that, in such a Poisson process, if we condition on the fact that there is at least one occurrence, then the largest event occurs at $1/(e-1)$ on average.

\subsection*{Acknowledgement} The authors thank Peter Winkler for useful comments.

\end{document}